\documentstyle[11pt]{article}
 \newcommand{\qb}[2]{{\left [{#1 \atop #2} \right]}}
 \newlength{\standardunitlength}
\newtheorem{cor}{Corollary} \newtheorem{lemma}{Lemma}
\newtheorem{theorem}{Theorem} \newtheorem{prop}{Proposition}
\newenvironment{proof}{\noindent {\sc Proof:}}{$\Box$ \vspace{2 ex}}

\begin{document}

\begin{center} Affine shuffles, shuffles with cuts, the Whitehouse
module, and patience sorting \end{center}

\begin{center}
By Jason Fulman
\end{center}

\begin{center}
Stanford University
\end{center}

\begin{center}
Department of Mathematics
\end{center}

\begin{center}
Building 380, MC 2125
\end{center}

\begin{center}
Stanford, CA 94305
\end{center}

\begin{center}
email:fulman@math.stanford.edu
\end{center}

\begin{center}
http://math.stanford.edu/$\sim$fulman
\end{center}

\begin{center}
1991 AMS Primary Subject Classifications: 20P05, 05E99
\end{center}

\newpage
Proposed running head: Affine shuffles

\newpage \begin{abstract} Using representation theoretic work on the
Whitehouse module, a formula is obtained for the cycle structure of a
riffle shuffle followed by a cut. It is proved that the use of cuts
does not speed up the convergence rate of riffle shuffles to
randomness. Type $A$ affine shuffles are compared with riffle shuffles
followed by a cut. Although these probability measures on the
symmetric group $S_n$ are different, they both satisfy a convolution
property. Strong evidence is given that when the underlying parameter
$q$ satisfies $gcd(n,q-1)=1$, the induced measures on conjugacy
classes of the symmetric group coincide. This gives rise to
interesting combinatorics concerning the modular equidistribution by
major index of permutations in a given conjugacy class and with a
given number of cyclic descents. Generating functions for the first
pile size in patience sorting from decks with repeated values are
derived. This relates to random matrices. \end{abstract}


\section{Introduction}

	In an effort to study the way real people shuffle cards, Bayer
and Diaconis \cite{BD} performed a definitive analysis of the
Gilbert-Shannon-Reeds model of riffle shuffling. For an integer $k
\geq 1$, a $k$-shuffle can be described as follows. Given a deck of
$n$ cards, one cuts it into $k$ piles with probability of pile sizes
$j_1,\cdots,j_k$ given by $\frac{{n \choose
j_1,\cdots,j_k}}{k^n}$. Then cards are dropped from the packets with
probability proportional to the pile size at a given time (thus if the
current pile sizes are $A_1,\cdots,A_k$, the next card is dropped from
pile $i$ with probability $\frac{A_i}{A_1+\cdots+A_k}$). It was proved
in \cite{BD} that $\frac{3}{2}log_2n$ shuffles are necessary and
suffice for a $2$-shuffle to achieve randomness (the paper \cite{A}
had established this result asymptotically in $n$). It was proved in
\cite{DMP} that if $k=q$ is a prime power, then the chance that a
permutation distributed as a $q$-shuffle has $n_i$ $i$-cycles is equal
to the probability that a uniformly chosen monic degree $n$ polynomial
over the field $F_q$ factors into $n_i$ irreducible polynomials of
degree $i$.

	These results have recently been extended to other Coxeter
groups and placed in a Lie theoretic setting. The paper \cite{BB}
defines hyperoctaheral shuffles using descent algebras and the paper
\cite{BHR} relates Gilbert-Shannon-Reeds shuffles to hyperplane
walks. The paper \cite{F2} defines riffle shuffling for arbitrary real
hyperplane arrangements, with convergence rates obtainable from the
theory in \cite{BHR}. The results of \cite{DMP} are given a Lie
theoretic formulation and extension, at least for types $A$ and $B$,
in \cite{F3} and \cite{F4}. (Random polynomials are replaced by the
semsimple orbits of the adjoint action of a finite group of Lie type
on its Lie algebra, and even in type $A$ restrictions on the
characteristic are needed). The paper \cite{F1} considers cycle
structure of permutations after biased shuffles, and generalizations
based on dynamical systems appear in \cite{La1},\cite{La2}.

	It is worth commenting that the combinatorics of type $A$
riffle shuffles is intimitely connected to cyclic and Hochschild
homology \cite{H},\cite{GS} to the Poincar\'e-Birkhoff-Witt theorem
\cite{BW}, free Lie algebras \cite{Ga}, and to Hopf algebras (Section
3.8 of \cite{SS}). In recent work, Stanley \cite{Sta} has related
biased riffle shuffles with representation theory of the symmetric
group, thereby giving an elementary probabilistic interpretation of
Schur functions and a different approach to some work in the random
matrix community. He recasts and extends results of \cite{BD} and
\cite{F1} using quasisymmetric functions.

	Using a construction of Cellini \cite{Ce1},\cite{Ce2}, the
paper \cite{F5} studies combinatorially much more intricate shuffles
called affine shuffles (they are reviewed in Section
\ref{preliminaries1}). The conjectures of \cite{F5} state in type $A$
that the chance that a permutation distributed as an affine $q$
shuffle has $n_i$ $i$-cycles is equal to the probability that a
uniformly chosen monic degree $n$ polynomial {\it with constant term
1} over the field $F_q$ factors into $n_i$ irreducible polynomials of
degree $i$ (the abstraction of these polynomials is semisimple
conjugacy classes of finite groups of Lie type). These conjectures are
remarkable in the sense that (unlike the Lie algebra case \cite{F3}),
no restrictions on the characteristic are needed and there seems to be
a reasonably natural way of associating to such a polynomial a
permutation in the right conjugacy class, such that choosing the
polynomial at random induces the affine shuffling measure. As emerges
in \cite{F5} (which gives an application to dynamical systems and
hints at number theoretic applications), this conjecture seems
challenging.

	The second type of shuffle to be studied in this paper is
riffle shuffling followed by a cut at a uniformly chosen random
position. Section \ref{preliminaries2} develops combinatorial
preliminaries of shuffles followed by cuts. It is shown there that
doing $r$ ``$k$-shuffles followed by a cut'' is the same as doing $r$
$k$-shuffles and then a single cut (this is known for $k=2$ from
\cite{Ce3}). It is proved that the total variation distance between a
sequence of $x$ riffle shuffles and $y$ cuts (performed in any order)
and the uniform distribution on $S_n$ is at least the total variation
distance between a sequence of $x$ riffle shuffles on $S_{n-1}$ and
the uniform distribution on $S_{n-1}$. In this precise sense cuts do
not help speed up riffle shuffles. This perhaps surprising fact can be
contrasted with a result of Diaconis \cite{D2}, who used
representation theory to show that although shuffling by doing random
tranpositions gets random in $\frac{1}{2}n log(n)$ steps, the use of
cuts at each stage drops the convergence time to $\frac{3}{8} n
log(n)$ steps. It would be worthwhile and interesting to
systematically study the effect of cuts on convergence rates of
shuffling methods. Section \ref{preliminaries2} also shows that a
riffle shuffle followed by a cut is at least as random (and sometimes
moreso) than a cut followed by a riffle shuffle.

	 Section \ref{representation} uses representation theoretic
work on the Whitehouse module to obtain a formula for the cycle
structure of a riffle shuffle followed by a cut. More precisely, it is
shown that the chance that a permutation distributed as a $q$ riffle
shuffle followed by a cut has $n_i$ $i$-cycles is equal to the
probability that a uniformly chosen monic degree $n$ polynomial with
{\it non-zero constant term} over the field $F_q$ factors into $n_i$
irreducible polynomials of degree $i$. The connection with
representation theory is illuminating; it follows for instance that
the work of \cite{DMP} on the cycle structure of riffle shuffles is
equivalent to a representation theoretic result of Hanlon
\cite{H}. The Whitehouse module is interesting in its own right and
appears in many places in mathematics; the interested reader should
see the transparencies for a talk on the Whitehouse module on Richard
Stanley's MIT website.

	Section \ref{compare} gives strong evidence for the assertion
that affine shuffles and shuffles followed by a cut, though different
probability measures, coincide when lumped according to conjugacy
classes provided that the prime power $q$ satisfies
$gcd(n,q-1)=1$. This leads to fascinating combinatorics concerning the
modular equidistribution by major index of permutations in a given
conjugacy class and with a given number of cyclic descents.
 
	Section \ref{patience} considers questions which turn out to
be related to cycle structure in multiset permutations. To motivate
things, we follow the recent preprint \cite{AD} in its description of
patience sorting (which can also be viewed as a toy model for
Solitaire, a card game which unlike BlackJack has been extremely
difficult to analyze mathematically). The simplest case is that one
starts with a deck of cards labelled $1,2,\cdots,n$ in random
order. Cards are turned up one at a time and dealt into piles on the
table according to the rule that a card is placed on the leftmost pile
whose top card is of higher value. If no such pile exists, the card
starts a new pile to the right. For example the ordering

\[ 7 5 1 3 6 2 4 \] leads to the arrangement

\[  \begin{array}{c c c c c}
		1 &  &   \\
		5 & 2 & 4  \\
		7 & 3 & 6
	  \end{array} \] The survey paper \cite{AD} details
connections of patience sorting with ideas ranging from random
matrices and the Robinson-Schensted correspondence to interacting
particle systems (\cite{AD2},\cite{BDJ} \cite{H},\cite{R}). For
example the number of piles in patience sorting is the length of the
longest increasing subsequence of a random permutation. After an
appropriate renormalization, this statistic has the same distribution
as the largest eigenvalue of a random $GUE$ matrix. The usefulness of
\cite{AD} is its insight that other functions of the random shape
obtained through patience sorting have interesting structure. They
give results for various pile sizes and suggest the search for
analogous results for the following two variations of patience sorting
from decks with repeated values: ties allowed (i.e. a card may be
placed on a card of the same value) or ties forbidden. Section
\ref{patience} gives generating functions for the first pile size for
these two variations, and also for the case where each card is chosen
at random from a finite alphabet. We remark that our approach also
works for a game one could call $m$-Solitaire in which each card value
may be placed on a given pile up to $m$ times (ties allowed and ties
forbidden corresponding to $m=\infty,m=1$ respectively). As these
games are of less interest and the calculation is grungier, we omit
it.

	It may strike the reader that Section \ref{patience}, though
involving cards and cycle structure, is unrelated to card
shuffling. However, there is a relation. To elaborate, suppose one
picks a random word of length $n$ from a totally ordered alphabet,
where the probability of getting symbol $i$ is $x_i$. The number of
piles in patience sorting applied to the word is the length of the
longest weakly increasing subsequence. Each word corresponds to a
possible riffle shuffle and the longest increasing subsequence in the
corresponding permutation is the longest weakly increasing subsequence
in the word \cite{Sta}.

	To close the introduction, we introduce some terminology that
will be used throughout the paper. $C_r(m)$ will denote the Ramanujan
sum $\sum_l e^{2\pi i l m/r}$ where the sum is over all $l$ less than
and prime to $r$. An element $w$ in the symmetric group $S_n$ is said
to have a descent at position $i$ if $1 \leq i \leq n-1$ and
$w(i)>w(i+1)$. The notation $d(w)$ will denote the number of descents
of $w$. The major index of $w$, denote $maj(w)$ will be the sum of the
positions $i$ at which $w$ has a descent. The permutation $w$ is said
to have a cyclic descent at $n$ if $w(n)>w(1)$. Then $cd(w)$ is
defined as $d(w)$ if $w$ has no cyclic descent at $n$ and as $d(w)+1$
if $w$ has a cyclic descent at $n$. As noted in \cite{Ce1}, the
concepts of descents and cyclic descents have analogs for all Weyl
groups. (The descents are the simple positive roots mapped to negative
roots and $w$ has a cyclic descent if it keeps the highest root
positive).

\section{Type $A$ affine shuffles}
\label{preliminaries1}

	This section gives six definitions for what we call type $A$
affine $k$-shuffles. The first two are due to Cellini \cite{Ce1} and
are the best in that they generalize to all Weyl groups. The next four
are due to the author \cite{F5} and are very useful for computational
purposes. (The final definition is a ``physical'' description of
affine type $A$ $2$-shuffles. A good problem is to extend the physical
description to higher values of $k$). In all of the definitions $k$ is
a positive integer. A seventh definition (which is conjectural but
very conceptual and potentially valid for all Weyl groups) appears in
the final section of \cite{F5}. As it requires some effort to describe
and will not be used here, we omit it.

{\bf Definitions of affine type $A$ $k$-shuffles}

\begin{enumerate}

\item Let $W_k$ be the subgroup of the type $A$
affine Weyl group generated by reflections in the $n$ hyperplanes
$\{x_1=x_2,\cdots,x_{n-1}=x_n,x_n-x_1=k\}$. This subgroup has index
$k^{n-1}$ in the affine Weyl group and has $k^{n-1}$ unique minimal
length coset representatives, each of which can be written as the
product of a permutation and a translation. Choose one of these
$k^{n-1}$ coset representatives uniformly at random and consider its permutation part. Define a type $A$ affine $k$-shuffle
to be the distribution on permutations so obtained.

\item A type $A$ affine $k$-shuffle assigns probability
to $w^{-1} \in S_n$ equal to $\frac{1}{k^{n-1}}$ multiplied by the number
of integers vectors $(v_1,\cdots,v_n)$ satisfying the conditions

\begin{enumerate}
\item $v_1+\cdots+v_n=0$
\item $v_1 \geq v_2 \geq \cdots \geq v_n, v_1-v_n \leq k$
\item $v_i>v_{i+1}$ if $w(i)>w(i+1)$ (with $1 \leq i \leq n-1$)
\item $v_1<v_n+k$ if $w(n)>w(1)$
\end{enumerate}

\item A type $A$ affine $k$-shuffle assigns probability
to $w^{-1} \in S_n$ equal to $\frac{1}{k^{n-1}}$ multiplied by the number
of partitions with $\leq n-1$ parts of size at most $k-cd(w)$ such
that the total number being partitioned has size congruent to $-maj(w)
\ mod \ n$. Equivalently, it assigns probability equal to
$\frac{1}{k^{n-1}}$ multiplied by the number of partitions with $\leq
k-cd(w)$ parts of size at most $n-1$ such that the total number being
partitioned has size congruent to $-maj(w) \ mod \ n$.

\item A type $A$ affine $k$-shuffle assigns probability
to $w^{-1} \in S_n$ equal to

\[ \begin{array}{ll} \frac{1}{n k^{n-1}} \sum_{r|n,k-cd(w)} {\frac{n+k-cd(w)-r}{r} \choose
\frac{k-cd(w)}{r}} C_r(-maj(w)) & \mbox{if \ $k-cd(w)>0$}\\
\frac{1}{k^{n-1}}  & \mbox{if $k-cd(w)=0, maj(w)=0$ mod n}\\
0 & \mbox{otherwise}.
\end{array} \]

\item A type $A$ affine $k$-shuffle assigns probability
to $w^{-1} \in S_n$ equal to

\[ \frac{1}{k^{n-1}} \sum_{r=0}^{\infty} Coeff. \ of \ q^{r \cdot n} \
in \ \left(q^{maj(w)} \qb{k+n-cd(w)-1}{n-1}\right),\] whre $\qb{A}{B}$
denotes the $q$-binomial coefficient
$\frac{(1-q)\cdots(1-q^A)}{(1-q)\cdots(1-q^B)
(1-q)\cdots(1-q^{A-B})}$.

\item A type $A$ affine $2$ shuffle has the following physical
description for the symmetric group $S_{2n}$.

	Step 1: Choose an even number between $1$ and $2n$ with the
probability of getting $2j$ equal to $\frac{{2n \choose
2j}}{2^{2n-1}}$. From the stack of $2n$ cards, form a second pile of
size $2j$ by removing the top $j$ cards of the stack, and then putting
the bottom $j$ cards of the first stack on top of them.

	Step 2: Now one has a stack of size $2n-2j$ and a stack of
size $2j$. Drop cards repeatedly according to the rule that if stacks
$1,2$ have sizes $A,B$ at some time, then the next card comes from
stack $1$ with probability $\frac{A}{A+B}$ and from stack 2 with
probability $\frac{B}{A+B}$. (This is equivalent to choosing uniformly
at random one of the ${2n \choose 2j}$ interleavings preserving the
relative orders of the cards in each stack).

	 The description of $x_2$ is the same for the symmetric
group $S_{2n+1}$, except that at the beginning of Step 1, the chance
of getting $2j$ is $\frac{{2n+1 \choose 2j}}{2^{2n}}$ and at the
beginning of Step 2, one has a stack of size $2n+1-2j$ and a stack of
size $2j$.

\end{enumerate}

	An important property of these shuffles is the so called
``convolution property'', which says that a $k_1$ shuffle followed by a
$k_2$ shuffle is equivalent to a $k_1k_2$ shuffle. It is interesting
that the type $A$ riffle shuffles of \cite{BD} satisfy the same
convolution property, as do some of the generalizations in \cite{F2}.

	The conjecture of \cite{F5} that relates the cycle structure
of permutations distributed as affine $q$-shuffles to the
factorization of polynomials with constant term 1 appears to be
interesting. For instance it is shown there that for the case of the
identity conjugacy class of $S_n$, it amounts to the $m=0$ case of
following observation in ``modular combinatorial reciprocity''. We
recently noticed that this reciprocity statement appears in an
invariant theoretic setting as Hermite reciprocity in \cite{EJ}.

{\it For any positive integers $x,y$, the
number of ways (disregarding order and allowing repetition) of writing
$m$ (mod $y$) as the sum of $x$ integers of the set $0,1,\cdots,y-1$
is equal to the number of ways (disregarding order and allowing
repetition) of writing $m$ (mod $x$) as the sum of $y$ integers of the
set $0,1,\cdots,x-1$} 

	Let $f_{n,k,d}$ be the coefficient of $z^n$ in
$(\frac{z^k-1}{z-1})^d$ and let $\mu$ be the Moebius function. Let
$n_i(w)$ be the number of $i$-cycles in a permutation $w$. Then
(loc. cit.) the conjecture is equivalent to the truly bizarre
assertion (which we intentionally do not simplify) that for all
$n,k$,

\begin{eqnarray*}
& & \sum_{m=0 \ mod \ n} Coef. \ of \ q^m u^n t^k \
in \ \sum_{n=0}^{\infty} \frac{u^n} {(1-tq)\cdots(1-tq^n)} \sum_{w
\in S_n} t^{cd(w)} q^{maj(w)} \prod x_i^{n_i(w)}\\ & = &
\sum_{m=0 \ mod \ k-1} Coef. \ of \ q^m u^n t^k \ in \
\sum_{k=0}^{\infty} t^k \prod_{i=1}^{\infty} \prod_{m=1}^{\infty}
(\frac{1}{1-q^mx_iu^i})^{1/i \sum_{d|i} \mu(d) f_{m,k,i/d}}.
\end{eqnarray*}

\section{Shuffles followed by a cut}
\label{preliminaries2}

	To begin we remark that although \cite{BD} offers a formula
for a shuffle followed by a cut, the formula is really for a cut
followed by a shuffle, which is different.

	Let $s$ be the element of the group algebra of $S_n$ denoting
a $k$-riffle shuffle. Let $\zeta$ be the cyclic permutation $(1 \cdots
n)$ and let $c=\frac{1}{n} \sum_{i=0}^{n-1} \zeta^n$. Thus in this
notation a shuffle followed by a cut is simply $cs$. The inverse of an
element $\sum r_w w$ of the group algebra will be taken to mean $\sum
r_w w^{-1}$.

	It is useful to recall the following formula of Bayer and
Diaconis.

\begin{theorem} \label{fromBD} (\cite{BD}) The coefficient of a
permutation $w$ in the element $s$ is \[ \frac{1}{k^n}
{n+k-d(w^{-1})-1 \choose n}.\] \end{theorem}

	Theorem \ref{correct} derives an analogous formula for a
shuffle followed by a cut.

\begin{theorem} \label{correct} The coefficient of a permutation $w$
in the element $cs$ is \[ \frac{1}{nk^{n-1}} {n+k-cd(w^{-1})-1 \choose
n-1}.\] \end{theorem}

\begin{proof} Consider instead the coefficient of $w$ in
$s^{-1}c$. This coefficient is equal to

\[ \frac{1}{n} \sum_{k=0}^{n-1} Coeff. \ of \ w \zeta^k \ in \ 
s^{-1}.\] The element $w \zeta^k$ maps $i$ to $w(i+k \ mod \ 
n)$. Consequently letting $cd(w)$ be the number of cyclic descents of
$w$, there are $cd(w)$ values of $k$ for which $w \zeta^k$ has
$cd(w)-1$ descents, and $n-cd(w)$ values of $k$ for which $w \zeta^k$
has $cd(w)$ descents. Combining this with Theorem \ref{fromBD} shows
that the coefficient of $w$ in $s^{-1}c$ is

\[ \frac{1}{nk^n} \left(cd(w) {n+k-cd(w) \choose n} + (n-cd(w))
{n+k-cd(w)-1 \choose n} \right),\] which simplifies to the formula in the
statement of the theorem. \end{proof}

	This yields the following combinatorial corollary.

\begin{cor} \label{countem} Let $B_{n,i}$ be the number of elements of
$S_n$ with $i$ cyclic descents. Let $A_{n,i}$ be the number of
elements of $S_n$ with $i-1$ descents. Then

\begin{enumerate}
\item $x^{n-1} = \sum_{1 \leq i \leq n-1} \frac{B_{n,i}}{n} {n+x-i-1
\choose n-1}$.

\item If $n > 1$ then $B_{n,i}=n A_{n-1,i}$.
\end{enumerate}
\end{cor} 

\begin{proof} The first assertion is immediate from Theorem
\ref{correct}. The second assertion follows from the first together
with the well-known facts that $A_{n,i}=A_{n,n+1-i}$ and that
$A_{n,i}$ is the unique sequence satisfying Worpitzky's identity $x^n
= \sum_{1 \leq i \leq n} A_{n,i} {x+i-1 \choose n}$. \end{proof}

	Theorem \ref{canmove} appears in \cite{Ce3} for the case
$k=2$. As noted there, it implies that $(cs)^h=cs^h$ for any natural
number $h$. The proof given here is simpler.

\begin{theorem} \label{canmove} $csc=cs$.
\end{theorem}

\begin{proof} Taking inverses and using the fact that $c^{-1}=c$, it
is enough to show that $cs^{-1}c=s^{-1}c$. The coefficient of $w$ in
$cs^{-1}c$ is

\[ \frac{1}{n} \sum_{k=0}^{n-1} Coeff. \ of \ \zeta^kw \ in \
s^{-1}c.\] It is easy to see that $cd(\zeta^kw)=cd(w)$ for all
$k$. The result now follows from Theorem \ref{correct}. \end{proof}

	Next recall the notion of total variation distance
$||P_1-P_2||$ between two probability distributions $P_1$ and $P_2$ on
a finite set $X$. It is defined as \[ \frac{1}{2} \sum_{x \in X}
|P_1(x)-P_2(x)|.\] The book \cite{D2} explains why this is a natural
and useful notion of distance between probability
distributions. $P_1*P_2$ (the convolution) is defined by
$P_1*P_2(\pi)=\sum_{\tau \in S_n} P_1(\pi \tau^{-1}) P_2(\tau)$, and
$P_1*\cdots*P_k$ is defined inductively. The following elementary (and
well known) lemma will be helpful.

\begin{lemma} \label{inequal} Let $P,Q$ be any measures on a finite
group $G$ and let $U$ be the uniform distribution on $G$. Then
$||P*Q-U|| \leq ||Q-U||$. \end{lemma}

	Theorem \ref{noneedcuts} shows that cuts do not speed up the
convergence rate of riffle shuffles.

\begin{theorem} \label{noneedcuts}
\begin{enumerate}

\item Let $S^{(k)},C,U$ denote the probability distribution
corresponding to a $k$-riffle shuffle, a cut, and the uniform
distribution respectively. Then $||C*S^{(k)}-U|| \leq ||S^{(k)}*C-U||$
and the inequality can be strict. (In words, a shuffle followed by a
cut is more random than a cut followed by a shuffle).

\item For $n>1$, $||C*S^{(k)}-U||_{S_n}=||S^{(k)}-U||_{S_{n-1}}$.

\item Let $W$ be the convolution of any finite sequence of riffle
shuffles and cuts. Let $W'$ be the convolution of the same finite
sequence, but with the cuts eliminated. (By abuse of notation, these
can be viewed on any symmetric group). Then

\[ ||W-U||_{S_n} \geq ||W'-U||_{S_{n-1}}. \]

\end{enumerate}
\end{theorem}

\begin{proof} For the first assertion, observe that Theorem
\ref{canmove} gives that $C*S^{(k)}=C*S^{(k)}*C$. Now use Lemma
\ref{inequal}. Computations with the symmetric group
$S_4$ show that the inequality can be strict.

	For the second assertion, let $B(n,i)$ be the number of
elements of $S_n$ with $i$ cyclic descents and let $A(n,i)$ be the
number of elements of $S_n$ with $i-1$ descents. Observe that

\begin{eqnarray*}
||C*S^{(k)}-U|| & = & \frac{1}{2} \sum_{i=1}^{n-1}  B(n,i) \left| \frac{{k+n-i-1 \choose n-1}}{nk^{n-1}}-\frac{1}{n!} \right|\\
& = & \frac{1}{2} \sum_{i=1}^{n-1} A(n-1,i) \left|\frac{{k+n-i-1 \choose n-1}}{k^{n-1}}-\frac{1}{(n-1)!} \right|\\
& = & ||S^{(k)}-U||_{S_n-1}.
\end{eqnarray*} The second equality is the second part of Corollary \ref{countem} and the final equality follows from Theorem \ref{fromBD}.

	For the third assertion, the inequality is clear if $W$ has no
cuts. Otherwise, combining the fact that $S^{(i)}*S^{(j)}=S^{(ij)}$
for any $i,j$ with Theorem \ref{canmove} shows that $W$ is equivalent
to a convolution of the form $S^{(k_1)}*C*S^{(k_2)}$ (with $k_1$ or
$k_2$ possibly $0$ and $S^{(0)}$ denoting the measure placing all mass
the identity). Now observe that

\begin{eqnarray*}
||S^{(k_1)}*C*S^{(k_2)}-U||_{S_n} & \geq & ||C*S^{(k_1)}*C*S^{(k_2)}-U||_{S_n}\\
& = & ||C*S^{(k_1k_2)}-U||_{S_n}\\
& = & ||S^{(k_1k_2)}-U||_{S_{n-1}}.
\end{eqnarray*} The first equality is Lemma \ref{inequal}, the second equality comes from Theorem \ref{canmove}, and the third equality is the second part of this theorem.
\end{proof}

	A formula for a cut followed by a riffle shuffle appears in
\cite{BD}, though it is not evident how it could be used to prove part
1 of Theorem \ref{noneedcuts}.

	As a final problem, we observe that the $n$-cycle
$\zeta=(1 \cdots n)$ is a minimal length Coxeter element for type
$A$. As there are analogs of shuffling for other finite Coxeter groups
\cite{BB},\cite{F2}, it may be possible to extend the results of this
paper to other Coxeter groups.

\section{Representation theory} \label{representation}

	This section uses representation theory to obtain a formula
for the cycle structure of a riffle shuffle followed by a cut.

	It is useful to recall the notion of a cycle index associated
to a character of the symmetric group. Letting $n_i(w)$ be the number
of $i$-cycles of a permutation $w$ and $N$ be a subgroup of $S_n$, one
defines $Z_N(\chi)$ as

\[ Z_N(\chi) = \frac{1}{|N|} \sum_{w \in N} \chi(w) \prod_i
a_i^{n_i(w)}.\] The cycle index stores complete information about the
character $\chi$. For a proof of the following attractive property of
cycle indices, see \cite{Fe}.

\begin{lemma} \label{Walter} Let $N$ be a subgroup of $S_n$ and $\chi$
a class function on $N$. Then 

\[ Z_{S_n}(Ind_N^{S_n}(\chi)) = Z_N(\chi).\]
\end{lemma}

	Next, recall that an idempotent $e$ of the group algebra of a
finite group $G$ defines a character $\chi$ for the action of $G$ on
the left ideal $KGe$ of the group algebra of $G$ over a field $K$ of
characteristic zero. For a proof of Lemma \ref{bridge}, which will
serve as a bridge between representation theory and computing measures
over conjugacy classes, see \cite{H}. For its statement, let $e<w>$ be
the coefficient of $w$ in the idempotent $e$.

\begin{lemma} \label{bridge} Let $C$ be a conjugacy class of the
finite group $G$, and let $\chi$ be the character associated to the idempotent
$e$. Then

\[ \frac{1}{|G|} \sum_{w \in C} \chi(w) = \sum_{w \in C} e<w>.\]
\end{lemma}

	 It is also convenient to define

\[ Z_{S_n}(e) = \sum_{w \in S_n} e<w> \prod_i a_i^{n_i(w)},\] which
makes sense for any element $e$ of the group algebra. Note that one
does not divide by the order of the group. When $e$ is idempotent and
$\chi$ is the associated character, Lemma \ref{bridge} can be
rephrased as

\[ Z_{S_n}(\chi) = Z_{S_n}(e).\]

	To proceed recall the Eulerian idempotents $e_n^j$,
$j=1,\cdots,n$ in the group algebra $QS_n$ of the symmetric group over
the rationals. These can be defined \cite{GS} as follows. Let
$s_{i,n-i}=\sum w$ where the sum is over all ${n \choose i}$
permutations $w$ such that $w(1)<\cdots<w(i)$, $w(i+1)<\cdots<w(n)$
and let $s_n=\sum_{i=1}^{n-1} s_{i,n-i}$. Letting $\mu_j=2^j-2$, the
$e_n^j$ are defined as

\[ e_n^j = \prod_{i \neq j} \frac{s_n-\mu_i}{(\mu_j-\mu_i)}.\] They
are orthogonal idempotents which sum to the identity.

	The following result, which we shall need, is due to
Hanlon. The symbol $\mu$ denotes the Moebius function of elementary
number theory.

\begin{theorem} (\cite{H}) \label{Hanl}

\[ 1+\sum_{n=1}^{\infty} \sum_{i=1}^n k^i Z_{S_n}(e_n^i) = \prod_{i
\geq 1} (1-a_i)^{-(1/i) \sum_{d|i} \mu(d) k^{i/d}}.\] \end{theorem}

\begin{theorem} (\cite{Ga}) \label{Gars}

\[ \sum_{i=1}^n k^i e_n^i = \sum_{w \in S_n} {n+k-d(w)-1 \choose n}
w.\] \end{theorem}

{\bf Remark:} Combining Lemma \ref{bridge} and Theorem \ref{Gars}, one
sees that the formula for the cycle structure of a riffle shuffle
\cite{DMP} and Theorem \ref{Hanl} imply each other. It is interesting
that both proofs used a bijection of Gessel and Reutenauer \cite{G}.

	To continue, we let $\overline{e_n^j}$ denote the idempotent
obtained by multiplying the coefficient of $w$ in $e_n^j$ by
$sgn(w)$. Let $\lambda_{n+1}$ be the $n+1$ cycle $(1 \ 2 \cdots n+1)$
and $\Lambda_{n+1}=\frac{1}{n+1} \sum_{i=0}^n (sgn \lambda_{n+1}^i)
\lambda_{n+1}^i$. Viewing $\overline{e_n^j}$ as in the group algebra
of $S_{n+1}$, Whitehouse \cite{W} proves that for $j=1,\cdots,n$ the
element $\Lambda_{n+1} \overline{e_n^j}$ is an idempotent in the group
algebra $QS_{n+1}$, which we denote by $f_{n+1}^j$. Whitehouse's main
result is the following:

\begin{theorem} (\cite{W}) \label{white} Let $F_{n+1}^j,\overline{E_{n}^j}$
be the irreducible modules corresponding to the idempotents
$f_{n+1}^j$ and $\overline{e_n^j}$. Then

\[ F_{n+1}^j \oplus \bigoplus_{i=1}^j \overline{E_{n+1}^i} =
\bigoplus_{i=1}^j Ind_{S_n}^{S_{n+1}} \overline{E_{n}^i}.\] \end{theorem}

	As final preparation for the main result of this section, we
link the idempotent $\Lambda_{n+1} \overline{e_n^j}$ with riffle shuffles
followed by a cut.

\begin{lemma} \label{form} The coefficient of $w$ in $\sum_{j=1}^n k^j
\Lambda_{n+1} \overline{e_n^j}$ is $sgn(w) \frac{1}{n+1} {k+n-cd(w)
\choose n}$. \end{lemma}

\begin{proof} Given Theorem \ref{Gars}, this is an elementary
combinatorial verification. \end{proof}

	Theorem \ref{cyc} now derives the cycle structure of a
permutation distributed as a shuffle followed by a cut. So as to
simplify the generating functions, recall that $\sum_{d|i} \mu(d)$
vanishes unless $i=1$.

\begin{theorem} \label{cyc}
\begin{eqnarray*}
&& 1+\sum_{n \geq 1} \sum_{w \in S_{n+1}}
\frac{1}{(n+1)k^{n+1}} {n+k-cd(w) \choose n} \prod_i a_i^{n_i(w)}\\
& = & 1-\frac{1}{k-1} -\frac{a_1}{k} + \frac{1}{k-1} \prod_{i \geq 1}
(1-\frac{a_i}{k^i})^{-1/i \sum_{d|i} \mu(d) (k^{i/d}-1)}.
\end{eqnarray*} If $k=q$ is the size of a finite field, this says that the cycle type of a
permutation distributed as a shuffle followed by the cut has the same
law as the factorization type of a monic degree $n$ polynomial over
$F_q$ with non-vanishing constant term.  \end{theorem}

\begin{proof} Replacing $a_i$ by $a_i k^i (-1)^{i+1}$, it is enough to
show that

\begin{eqnarray*}
& & 1+\sum_{n \geq 1} \sum_{w \in S_{n+1}} sgn(w) \frac{1}{(n+1)}
{n+k-cd(w) \choose n} \prod_i a_i^{n_i(w)}\\
& = &1-\frac{1}{k-1} -a_1 +
\frac{1}{k-1} \prod_{i \geq 1} (1-(-1)^{i+1}a_i)^{-1/i \sum_{d|i}
\mu(d) (k^{i/d}-1)}.
\end{eqnarray*}

	Using Lemmas \ref{Walter}, \ref{bridge}, \ref{form} and
Theorem \ref{white}, one sees that

\begin{eqnarray*}
& &  1+\sum_{n \geq 1} \sum_{w \in S_{n+1}} sgn(w) \frac{1}{(n+1)}
{n+k-cd(w) \choose n} \prod_i a_i^{n_i(w)}\\
& =& 1+\sum_{n=1}^{\infty} \sum_{j=1}^n k^j Z_{S_{n+1}}(f_n^j)\\
& = & 1 + \sum_{n=1}^{\infty} \sum_{j=1}^n k^j \sum_{i=1}^j Z_{S_{n+1}}(Ind_{S_n}^{S_{n+1}}(\overline{e_n^i})) - \sum_{n=1}^{\infty} \sum_{j=1}^n k^j \sum_{i=1}^j Z_{S_{n+1}}(\overline{e_{n+1}^i})\\
& = & 1+a_1 \sum_{n=1}^{\infty} \sum_{i=1}^n Z_{S_n}(\overline{e_n^i})(\frac{k^{n+1}-k^i}{k-1}) - \sum_{n=1}^{\infty} \sum_{i=1}^n Z_{S_{n+1}}(\overline{e_{n+1}^i}) (\frac{k^{n+1}-k^i}{k-1})\\
& = &  1+a_1 k Z_{S_1}(\overline{e_1}) + \frac{a_1k-1}{k-1} \sum_{n=2}^{\infty} k^n \sum_{i=1}^n Z_{S_n}(\overline{e_n^i}) + \frac{1-a_1}{k-1} \sum_{n=2}^{\infty} \sum_{i=1}^n k^i Z_{S_n}(\overline{e_n^i}).
\end{eqnarray*} To simplify things further, recall that $\sum_{i=1}^n
Z_{S_n}(\overline{e_n^i})$ is $a_1^n$ since the $\overline{e_n^i}$'s sum to the
identity. The above then becomes

\[ 1-\frac{1}{k-1}-a_1+\frac{1-a_1}{k-1} (1+\sum_{n=1}^{\infty}
\sum_{i=1}^n k^i Z_{S_n}(\overline{e_n^i})),\] so the sought result follows
from Theorem \ref{Hanl}.
\end{proof}

	Before continuing, we observe that a combinatorial proof of
Theorem \ref{cyc} (which must exist) would give a new proof of Theorem
\ref{white}, by reversing the steps.

	Upon hearing about Theorem \ref{cyc}, Persi Diaconis
immediately asked for the expected number of fixed points after a
$k$-riffle shuffle followed by a cut, suggesting that it should be
smaller than for a $k$ riffle shuffle. Using the methods of Section 5
of \cite{DMP}, one can readily derive analogs of all of the results
there. As an illustrative example, Corollary \ref{fixpoint} shows that
the expected number of fixed points after a $k$-riffle shuffle
followed by a cut is the same as for a uniform permutation, namely 1
(the answer for $k$-riffle shuffles is $1+1/k+\cdots+1/k^{n-1}$). Two
other examples are worth mentioning and will be treated in Corollary
\ref{others}.

\begin{cor} \label{fixpoint} The expected number of fixed points after
a $k$-riffle shuffle followed by a cut is 1. \end{cor}

\begin{proof} The case $n=1$ is obvious. Multiplying $a_i$ by $u$ in
the statement of Theorem \ref{cyc} shows that

\begin{eqnarray*}
&& 1+\sum_{n \geq 1} \sum_{w \in S_{n+1}} u^{n+1}
\frac{1}{(n+1)k^{n+1}} {n+k-cd(w) \choose n} \prod_i a_i^{n_i(w)}\\
& = & 1-\frac{1}{k-1} -\frac{ua_1}{k} + \frac{1}{k-1} \prod_{i \geq 1}
(1-\frac{u^i a_i}{k^i})^{-1/i \sum_{d|i} \mu(d) (k^{i/d}-1)}.
\end{eqnarray*} To get the generating function in $u$ (for $n \neq 1$) for the expected number of fixed points in a riffle shuffle followed by a cut, one multiplies the right hand side by $k$, sets $a_2=a_3=\cdots=1$, differentiates with respect to $a_1$, and then sets $a_1=1$. Doing this yields the generating function \[ -u + u \prod_{i \geq 1} (1-\frac{u^i}{k^i})^{-1/i \sum_{d|i} \mu(d) k^{i/d}}.\] The result now follows from the identity

\[ \prod_{i \geq 1} (1-\frac{u^i}{k^i})^{-1/i \sum_{d|i} \mu(d)
k^{i/d}} = \frac{1}{1-u},\] which is equivalent to the assertion that
a monic degree $n$ polynomial over $F_q$ has a unique factorization
into irreducibles, since $1/i \sum_{d|i} \mu(d) k^{i/d}$ is the number
of irreducible polynomials of degree $i$ over the field $F_k$. \end{proof}

\begin{cor} \label{others} Fix $u$ with $0<u<1$. Let $N$ be chosen
from $\{0,1,2,\cdots\}$ according to the rule that $N=0$ with
probability $\frac{1-u}{1-u/k}$ and $N=n \geq 1$ with probability
$\frac{(k-1)(1-u)u^n}{k-u}$. Given $N$, let $w$ be the result of a
random $k$ shuffle followed by a cut. Let $N_i$ be the number of
cycles of $w$ of length $i$. Then the $N_i$ are independent and $N_i$
has a negative binomial distribution with parameters $1/i \sum_{d|i}
\mu(d) (k^{i/d}-1)$ and $(u/k)^i$. Consequently, for fixed $k$ as $n
\rightarrow \infty$, the joint distribution of the number of $i$
cycles after a $k$-shuffle followed by a cut converges to independent
negative binomials with parameters $1/i \sum_{d|i} \mu(d) (k^{i/d}-1)$
and $(1/k)^i$. \end{cor}

\begin{proof} Theorem \ref{cyc} and straightforward manipulations
give that

\begin{eqnarray*}
&& 1+ \frac{k-1}{k} \sum_{n \geq 1} \sum_{w \in S_{n}}
\frac{u^n}{nk^{n-1}} {n+k-cd(w)-1 \choose n-1} \prod_i a_i^{n_i(w)}\\
& = & \prod_{i \geq 1}
(1-\frac{a_i u^i}{k^i})^{-1/i \sum_{d|i} \mu(d) (k^{i/d}-1)}.
\end{eqnarray*} Setting all $a_i=1$ gives the equation

\[ 1+\frac{(k-1)u}{k(1-u)} = \prod_{i \geq 1}
(1-\frac{u^i}{k^i})^{-1/i \sum_{d|i} \mu(d) (k^{i/d}-1)}.\] Taking
reciprocals and multiplying by the first equation gives

\begin{eqnarray*}
&& (\frac{1-u}{1-u/k})+ \frac{(k-1)(1-u)}{k-u} \sum_{n \geq 1} \sum_{w \in
S_{n}} \frac{u^n}{nk^{n-1}} {n+k-cd(w)-1 \choose n-1} \prod_i
a_i^{n_i(w)}\\
& = & \prod_{i \geq 1} (\frac{1-\frac{u^i}{k^i}}{1-\frac{a_i u^i}{k^i}})^{1/i
\sum_{d|i} \mu(d) (k^{i/d}-1)}, \end{eqnarray*} proving the first
assertion of the corollary.

	For the second assertion there is a technique simpler than
that in \cite{DMP}. Rearranging the last equation gives that

\begin{eqnarray*}
&& (\frac{1-u}{1-1/k})+ \sum_{n \geq 1} \sum_{w \in
S_{n}} \frac{(1-u) u^n}{nk^{n-1}} {n+k-cd(w)-1 \choose n-1} \prod_i
a_i^{n_i(w)}\\
& = & \frac{1-u/k}{1-1/k} \prod_{i \geq 1} (\frac{1-\frac{u^i}{k^i}}{1-\frac{a_i u^i}{k^i}})^{1/i
\sum_{d|i} \mu(d) (k^{i/d}-1)}. \end{eqnarray*} Letting $g(u)$ be a generating function with a convergent Taylor series, the limit coefficient of $u^n$ in $\frac{g(u)}{1-u}$ is simply $g(1)$. This proves the second assertion. \end{proof}

\section{Conjugacy classes} \label{compare}

	The aim of this section is to give evidence for the conjecture
at the end of Section \ref{preliminaries1}, in the case when
$gcd(q-1,n)=1$. Note that under this assumption, a uniformly chosen
degree $n$ polynomial with non-zero constant term and a uniformly
chosen degree $n$ polynomial with constant term $1$ have the same
chance of factoring into $n_i$ $i$-cycles. Hence in this case the
conjecture amounts to the assertion that affine shuffles and shuffles
followed by a cut, though different probability measures, induce the
same distribution on conjugacy classes. Before posing a problem which
would explain why this should hold, some lemmas are needed.

\begin{lemma} \label{sum} If $r>1$, then $\sum_{j=0}^{r-1} C_r(j)=0$.
\end{lemma}

\begin{proof} If $l$ is relatively prime to $r$, then multiplication
by $l$ permutes the numbers $\{0,1,\cdots,r-1\}$ mod $r$. Thus
 \[ \sum_{j=0}^{r-1} C_r(j) = \sum_{0 \leq l \leq r \atop gcd(l,r)=1}
\sum_{j=0}^{r-1} e^{2\pi i j l/r} = \phi(r) \sum_{j=0}^{r-1} e^{2\pi i
j /r}=0.\] \end{proof}

\begin{lemma} \label{divide} For $n \geq 1$, let $t$ be the largest
divisor of $n$ such that $gcd(cd-1,t)=1$. Suppose that $gcd(n,q-1)=1$
and that $r$ divides $n$ and $q-cd$. Then $r$ divides $t$. \end{lemma}

\begin{proof} Observe that $gcd(r,cd-1)=1$. For suppose there is some
$a>1$ dividing $r$ and $cd-1$. Then $a$ divides $q-cd$ and $cd-1$,
hence $q-1$. Since $a$ divides $r$ and $r$ divides $n$, it follows
that $a$ divides $n$. This contradicts the assumption that
$gcd(q-1,n)=1$. \end{proof}

	Next we pose the problem of determining whether or not the
following statement holds. 

{\bf Statement 1:} For $n \geq 1$, let $t$ be the largest divisor of
$n$ such that $gcd(cd-1,t)=1$. Then for every conjugacy class $C$ of
$S_n$, the set of permutations in $C$ with $cd$ cyclic descents has
its major index equidistributed mod $t$.

	Theorem \ref{imply} shows that if Statement 1 holds, then the
conjecture about the cycle structure of permutations distributed as
affine shuffles is correct. Some evidence in favor of Statement 1 is
then given. 

\begin{theorem} \label{imply} Suppose that $gcd(q-1,n)=1$. If
Statement 1 is correct, then affine shuffles and shuffles followed by
a cut have the same distribution on conjugacy classes. \end{theorem}

\begin{proof} Suppose that Statement 1 is correct and recall the third
definition of affine $q$-shuffles in Section \ref{preliminaries1}. If
$q<cd(w)$ then both the affine $q$-shuffle and the $q$-riffle shuffle
followed by a cut assign probability $0$ to $w$. If $q=cd(w)$, then
the affine $q$-shuffle assigns probability $\frac{1}{q^{n-1}}$ to $w$
if $maj(w)=0 \ mod \ n$, and 0 otherwise. If $q=cd(w)$, then the
$q$-riffle shuffle followed by a cut associates probability
$\frac{1}{nq^{n-1}}$ to $w$. Since $q=cd$, the $t$ in Statement 1 is
equal to $n$, which implies that for every conjugacy class $C$, the
set of permutations in $C$ with $cd$ cyclic descents has major index
equidistributed mod $n$. Hence Statement 1 holds in this case.

	The third and final case is that $q>cd(w)$. Suppose that $r>1$
divides $n$ and $q-cd$. Lemma \ref{divide} implies that $r$ divides
$t$. Hence by Statement 1, for any conjugacy class $C$, the set of
permutations with $cd$ cyclic descents has its major index
equidistributed mod $r$. Consequently (the second equality below
holding by Lemma \ref{sum} and the equidistribution property mod $r$),

\begin{eqnarray*}
& &\sum_{cd=1}^n \sum_{w \in C \atop cd(w)=cd} \frac{1}{n q^{n-1}}
\sum_{r|n,q-cd} {\frac{n+q-cd(w)-r}{r} \choose \frac{q-cd(w)}{r}}
C_r(-maj(w))\\
& = & \sum_{cd=1}^n \frac{1}{n q^{n-1}} \sum_{r|n,q-cd} {\frac{n+q-cd-r}{r} \choose \frac{q-cd}{r}} \sum_{w \in C \atop cd(w)=cd} C_r(-maj(w))\\
& = & \sum_{cd=1}^n \frac{1}{n q^{n-1}} {n+q-cd-1 \choose n-1} \sum_{w \in C \atop cd(w)=cd} 1\\
& = & \sum_{w \in C} \frac{1}{n q^{n-1}} {n+q-cd(w)-1 \choose n-1}.
\end{eqnarray*} From Theorem \ref{correct} (the formula for a $q$-riffle shuffle followed by a cut), Statement 1 follows.
\end{proof} 

	Next we consider evidence in favor of the ideas of this
section. Incidentally, given Section \ref{representation}, Proposition
\ref{obvious} confirms Conjecture 1 of \cite{F5} (in type $A$) when
$n$ is prime and $q$ is a power of $n$.

\begin{prop} \label{obvious} Suppose that $n$ is prime and that $q$ is
a power of $n$. Then type $A$ affine $q$-shuffles are exactly the same
as $q$-riffle shuffles followed by a cut. \end{prop}

\begin{proof} The probability that an affine $q$ shuffle yields $w$ is
 \[ \frac{1}{n q^{n-1}} \sum_{r|n,q-cd(w^{-1})}
{\frac{n+q-cd(w^{-1})-r}{r} \choose \frac{q-cd(w^{-1})}{r}}
C_r(-maj(w^{-1}))\] Since $1 \leq cd(w) \leq n-1$ for any $w$ in
$S_n$, the assumptions on $n$ and $q$ imply that the only $r$ dividing
$n$ and $q-cd(w^{-1})$ is $r=1$. The result now follows from Theorem
\ref{correct}. \end{proof}

\begin{theorem} Statement 1 holds for the identity conjugacy class and
for the conjugacy class of simple transpositions. \end{theorem}

\begin{proof} For the identity conjugacy class, use the third
definition of affine $q$-shuffles in Section \ref{preliminaries1},
together with the assumption that $gcd(n,q-1)=1$.

	Next consider the case of simple transpositions. Suppose that
$n \geq 4$, the other cases being trivial. One checks that all simple
transpositions $(i,j)$ with $i<j$ have either $2$ or $3$ cyclic
descents. The easy case is that of $2$ cyclic descents. The possible
values of $(i,j)$ are then $(i,i+1)$ for $1 \leq i \leq n-2$, $(1,n)$
and $(n-1,n)$. The values of the major index so obtained are
$\{1,\cdots,n\}$ and each value is hit once. Thus Statement 1 holds in
this case.

	The harder case is that of 3 cyclic descents. The relevant
transpositions are $(i,j)$ with $1 \leq i,j<n$ and $j \neq i+1$
(having major index $i+j-1$) and $(i,n)$ with $2 \leq i<n-1$ (having
major index $i+n-1$).

	First suppose that $n$ is odd. It suffices to prove that
$\sum_{w=(i,j) \atop cd(w)=3} x^{maj(w) \ mod \ n}$ is a multiple of
$\frac{x^n-1}{x-1}$. Calculating gives

\begin{eqnarray*}
\sum_{w=(i,j) \atop cd(w)=3} x^{maj(w) \ mod \ n} & = & \sum_{i=1}^{(n-3)/2} x^i \sum_{j=i+1}^{n-i-1} x^j +
\frac{1}{x^n} \sum_{i=1}^{(n-3)/2} x^i \sum_{j=n-i}^{n-2} x^j\\
&&  + \frac{1}{x^n}
\sum_{i=(n-1)/2}^{n-3}x^i \sum_{j=i+1}^{n-2} x^j + \sum_{i=1}^{n-3} x^i\\
& = & \frac{1}{x-1}\left( \sum_{i=1}^{(n-3)/2} x^i (x^{n-i}-x^{i+1}) + \frac{1}{x^n} \sum_{i=1}^{(n-3)/2} x^i (x^{n-1}-x^{n-i}) \right)\\
& & + \frac{1}{x-1}\left(\frac{1}{x^n} \sum_{i=(n-1)/2}^{n-3} x^i (x^{n-1}-x^{i+1}) + x^{n-2}-x \right)\\
& = & \frac{n-3}{2} \frac{x^n-1}{x-1},
\end{eqnarray*} as desired.

	Next suppose that $n=2^a$ with $a>0$. It suffices to prove
that $\sum_{w=(i,j) \atop cd(w)=3} x^{maj(w) \ mod \ n}$ is a
polynomial multiple of $\frac{x^{n/2^a}-1}{x-1}$. Calculating as above
(and omitting the steps analogous to the previous computation) gives
that

\begin{eqnarray*}
\sum_{w=(i,j) \atop cd(w)=3} x^{maj(w) \ mod \ n} &=& \frac{1}{x-1} (1+x^2+x^4+\cdots+x^{n-2}-x-x^3-x^5-\cdots-x^{n-1})\\
&=& -\frac{x^n-1}{x^2-1}\\
&=& -\frac{x^{n/2^a}-1}{x-1} \frac{x^{(2^a-1)n/2^a}+\cdots+x^{n/2^a}+1}{x+1}.
\end{eqnarray*} Since $n/2^a$ is odd, it follows that $x^{(2^a-1)n/2^a}+\cdots+x^{n/2^a}+1$ is divisible by $x+1$.
\end{proof}

\section{Patience sorting} \label{patience}

	Having described the motivation in the introduction, we
outline and then execute a strategy for obtaining generating function
information for the first pile size in patience sorting from decks
with repeated values. The first step is to apply ideas of Foata to
obtain generating functions for multiset permutations by the number of
cycles. The second step is to give a multiset records-to-cycles
bijection (generalizing the one used in \cite{AD}), which converts
information about the distribution of cycles to information about the
distribution of records. The final step is to read information off of
the generating function.

	Some notation is needed. Let $\vec{a}$ denote the vector
$(a_1,a_2,\cdots)$ with $a_i \geq 0$ and $\sum a_i < \infty$. Let
$Mult(\vec{a})$ denote the collection of all ${\sum a_i \choose a_1,
a_2, \cdots}$ words of length $\sum a_i$ formed from $a_i$ $i$'s.

	We recall Foata's theory of cycle structure for multisets
\cite{Fo}, following Knuth's superb exposition \cite{Kn}. Suppose that
the elements of the multiset are linearly ordered. Then multiset
permutations can be written in two-line notation

\[ \left( \begin{array}{c c c c c c c c c c}
		a & a & a & b& b& c& d& d& d& d \\
		c & a & b & d& d& a& b& d& a& d
	  \end{array} \right) \] Foata introduced an intercalation
product $_T$ which multiplies two multiset permutations $\alpha$ and
$\beta$ by expressing $\alpha$ and $\beta$ in two line notation,
juxtaposing these two-line notations and then sorting the columns in
non-decreasing order of the top line. For example

\[  \begin{array}{c c c c c c c c c c c c c c c c c c c c c c c c}
		a & a & b & c& d& & & & a& b & d & d & d & _= & a&a&a&b&b&c&d&d&d&d \\
		c & a & d & a& b& & _T  & & b& d &d & a & d &  & c&a&b&d&d&a&b&d&a&d
	  \end{array}  \] Foata proved that if the elements of the multiset $M$ be
linearly ordered by the relation $<$, then the permutations $\pi$ of $M$
correspond exactly to the possible intercalations
\[ \pi = (x_{11} \cdots x_{1n_1} y_1) _T (x_{21} \cdots x_{2n_2} y_2)
_T \cdots _T (x_{t1} \cdots x_{tn_t} y_t) \] with $y_1 \leq y_2 \cdots
\leq y_t$ and $y_i <x_{ij}$ for $1 \leq j \leq n_i$, $1 \leq i \leq
t$. This defines a notion of cycle structure for multiset permutations by letting the cycles be the intercalation factors.  Let
$C_i(\pi)$ be the number of length $i$ cycles of $\pi$ and $C(\pi) = \sum C_i(\pi)$ be the total
number of cycles. Let $C_i'(\pi)$ be the number of $i$-cycles of
$\pi$, where cycles with the same minimum value $y_j$ are counted at
most once. Let $C'(\pi) = \sum C_i'(\pi)$. For example the multiset
permutation

\[ (4 3 1) _T (2 3 1) _T (4) \] satisfies $C_3(\pi)=2$,
$C_3'(\pi)=1$, $C(\pi)=3$, and $C'(\pi)=2$. 

	Proposition \ref{genfunc} gives generating functions for
multiset permutations.

\begin{prop} \label{genfunc}  \[ 1+\sum_{\vec{a}} \sum_{\pi \in
Mult(\vec{a})} u^{C(\pi)} \prod_{i \geq 1} x_i^{a_i} =
\prod_{k=1}^{\infty} \frac{1}{1-\frac{x_k u}{1-\sum_{j>k} x_j}} \]

\[ 1+\sum_{\vec{a}} \sum_{\pi \in Mult(\vec{a})} u^{C'(\pi)} \prod_{i
\geq 1} x_i^{a_i} = \prod_{k=1}^{\infty}
\left(1+u\frac{\frac{x_k}{1-\sum_{j>k} x_j}}{1-\frac{x_k}{1-\sum_{j>k}
x_k}} \right) \] \end{prop}

\begin{proof} Both generating functions follow easily from Foata's
method of representing permutations by intercalations. The $k$'s on
the right hand-side index the letters of the alphabet. The point is
that cycles are formed by fixing a smallest element $k$ and specifying
an ordered choice of elements larger than $k$; permutations are
ordered multisets of such cycles.\end{proof}

	We remark that the generating functions of Proposition
\ref{genfunc} are quasi-symmetric functions in the sense that for
any $i_1 < \cdots < i_n$ and $j_1 < \cdots < j_n$ the coefficients of
$x_1^{i_1} \cdots x_n^{i_n}$ and $x_1^{j_1} \cdots x_n^{j_n}$ are
equal.

	Theorem \ref{convert} converts information about the
distribution of cycles to information about the distribution of
records. Some further notation is needed for its statement. Let
$\pi^{rev}$ be the word obtained by reading from right to left the
bottom line in the $2$-line form of $\pi$. Recalling the definition of
Solitaire from the introduction, let $P_i(\pi)$ and $P_i'(\pi)$ be the
number of cards in pile $i$ of Solitaire with ties allowed and ties
forbidden respectively.

\begin{theorem} \label{convert}
\begin{enumerate}

\item For any given $\vec{a}$, there is a bijection $\Phi:
Mult(\vec{a}) \mapsto Mult(\vec{a})$ such that if $R_1,\cdots,R_t$ are
the positions of the left-to-right minima (ties are allowed) of
$\pi^{rev}$, then $R_2-R_1,R_3-R_2,\cdots,R_t-R_{t-1},(\sum
a_i)+1-R_t$ are the cycle lengths in Foata's factorization of
$\Phi(\pi)$.

\item $P_1(\pi) = C(\Phi(\pi))$ and $P_1'(\pi) = C'(\Phi(\pi))$.
\end{enumerate}
\end{theorem}

\begin{proof} For the first assertion, define $\Phi(\pi)$ as an
intercalation of cycles formed by entries in the bottom line of $\pi$,
with cycles (from left-to-right) having lengths $(\sum a_i)+1-R_t,
R_t-R_{t-1},\cdots,R_2-R_1$. The assertion is then evident, and the
following example may help to untangle the notation. The multiset
permutation $\pi = d \ d \ b \ c \ d \ b \ b \ c \ a\ b\ a\ c\ d\ b\
d$ has $\pi^{rev} = d\ b\ d\ c\ a\ b\ a\ c\ b\ b\ d\ c\ b\ d\ d$ with
$R_1(\pi^{rev})=1$, $R_2(\pi^{rev})=2$, $R_3(\pi^{rev})=5$,
$R_4(\pi^{rev})=7$, and $\sum a_i+1=16$. Thus forming from $\pi$
cycles of lengths $9,2,3,1$ gives $\Phi(\pi)$ as the intercalation

\[ (d\ d\ b\ c\ d\ b\ b\ c\ a)_T (b\ a)_T (c\ d\ b)_T (d).\]

	For the second assertion, we give the argument for the first
equality, the argument for the second assertion being analogous. The
point is that the number of cards in pile 1 of Solitaire with ties
allowed applied to $\pi$ is simply the number of left-to-right minima
with ties allowed of $\pi$. The result now follows from the first
assertion.
\end{proof}

	Proposition \ref{words} shows that when one considers random
words from a finite alphabet, there is a factorization for the full
cycle structure vector, not only the number of cycles. Recent work of
Tracy and Widom \cite{TW} connects random words with random matrices
chosen from the Laguerre ensemble. Let $Word_n(N)$ be the $N^n$ words
of length $n$ from an alphabet on $N$ letters (say
$1,2,\cdots,N$). Each such word can be viewed as a multiset
permutation.

\begin{prop} \label{words}
\[ 1+\sum_n \frac{1}{N^n} \sum_{\pi \in Word_n(N)} \prod_{k=1}^N
x_k^{a_k(\pi)}\prod_{i \geq 1} u_i^{C_i(\pi)} = \prod_{i \geq 1}
\prod_{k=1}^N \left(\frac{1}{1-u_i \frac{x_k}{N} (\sum_{j=k+1}^N
\frac{x_j}{N})^{i-1}}\right).\]

\[ (1-x)+\sum_n \frac{(1-x)x^n}{N^n} \sum_{\pi \in Word_n(N)}
\prod_{i \geq 1} u_i^{C_i(\pi)} = \prod_{i \geq 1}
\prod_{k=1}^N \frac{1-\frac{x^i}{N} (\frac{N-k}{N})^{i-1}}{1-\frac{u_i x^i}{N}
(\frac{N-k}{N})^{i-1}}.\]

\[ 1+\sum_n \frac{1}{N^n} \sum_{\pi \in Word_n(N)} \prod_{k=1}^N
x_k^{a_k(\pi)}\prod_{i \geq 1} u_i^{C_i'(\pi)} = \prod_{i \geq 1}
\prod_{k=1}^N \left( 1+u_i \frac{x_k}{N} \frac{1}{1-\sum_{j=k+1}^N \frac{x_j}{N}}\right).\]

\[ (1-x)+\sum_n \frac{(1-x)x^n}{N^n} \sum_{\pi \in Word_n(N)}
 \prod_{i \geq 1} u_i^{C_i'(\pi)} = \prod_{i \geq 1}
\prod_{k=1}^N \frac{1+u_i \frac{x}{N} \frac{1}{1-\sum_{j=k+1}^N
\frac{x}{N}}}{1+ \frac{x}{N} \frac{1}{1-\sum_{j=k+1}^N
\frac{x}{N}}}.\]
\end{prop}

\begin{proof} For the first assertion, note by Foata's representation
of multiset permutations as intercalations that each $i$-cycle is
formed by fixing a smallest element $k$ and specifying an ordered
choice of $i-1$ elements larger than $k$ to occupy the first $i-1$
positions of the cycle. Since multiset permutations are ordered
multisets of such cycles, one concludes that \[ 1+\sum_n \sum_{\pi \in
Word_n(N)} \prod_{k=1}^N x_k^{a_k(\pi)}\prod_{i \geq 1} u_i^{C_i(\pi)}
= \prod_{i \geq 1} \prod_{k=1}^N \frac{1}{1-u_i x_k (\sum_{j=k+1}^N
x_j)^{i-1}}.\] Now replace each $x_i$ by $\frac{x_i}{N}$.

	To prove the second assertion, replacing each $x_i$ by $x$ in
the first yields the equation \[ 1+\sum_n \frac{1}{N^n} \sum_{\pi \in
Word_n(N)} x^n \prod_{i \geq 1} u_i^{C_i(\pi)} = \prod_{i \geq 1}
\prod_{k=1}^N \frac{1}{1-\frac{u_i x^i}{N} (\frac{N-k}{N})^{i-1}}.\]
Setting all $u_i=1$ and taking reciprocals shows that \[ 1-x =
\prod_{i \geq 1} \prod_{k=1}^N \left(1- \frac{x^i}{N}
(\frac{N-k}{N})^{i-1}\right).\] The result follows by multiplying the
previous two equations.

	The arguments for the third and fourth assertions are
analogous.
\end{proof}

	The second and fourth equations have probabilistic
interpretations. For instance in the second equation, fix $x$ such
that $0<x<1$. The equation then says that if one picks $n$
geometrically with probability $(1-x)x^n$ and then picks $\pi \in
Word_n(N)$ uniformly at random, the random variables $C_i(\pi)$ are
sums of independent geometrics. In the fourth equation the $C_i'$
become sums of independent binomials. Lemma \ref{bign} permits
asymptotic statements in the $n \rightarrow \infty$ limit.

\begin{lemma} \label{bign} If $f(1)<\infty$ and the Taylor series of
$f$ around 0 converges at $u=1$, then

\[ lim_{n \rightarrow \infty} [u^n] \frac{f(u)}{1-u} = f(1). \]
\end{lemma}

\begin{proof}Write the Taylor expansion $f(u) = \sum_{n=0}^{\infty}
a_n u^n$. Then observe that $[u^n] \frac{f(u)}{1-u} = \sum_{i=0}^n
a_i$.  \end{proof}

	As a corollary, one sees for instance that as $n \rightarrow
\infty$, the number of $i$-cycles of a random length $n$ word from the
alphabet $\{1,\cdots,N\}$ converges to a sum of independent geometrics
with parameters $(1-\frac{k}{N})^{i-1}$ as $k=1,\cdots,N$. For more on
this type of factorization result and its applications, see \cite{LS}
for the symmetric groups, \cite{DS} for the compact classical groups,
and \cite{F6} for the finite classical groups.

	Finally, we consider the application of Proposition
\ref{genfunc} and Theorem \ref{convert} to patience sorting. As above,
$P_1(\pi)$ and $P_1'(\pi)$ be the number of cards in pile $1$ of
patience sorting with ties allowed and ties forbidden respectively.

\begin{theorem} Let $\pi$ be chosen uniformly at random from the
possible orderings of a deck of cards with $a_i$ cards labelled
$i$. Then $E(P_1) = \sum_{k} \frac{a_k}{a_1+\cdots+a_{k-1}+1}$ and
$E(P_1') = \sum_k \frac{a_k}{a_1+\cdots+a_{k-1}+a_k}$.  \end{theorem}

\begin{proof} By Proposition \ref{genfunc} and Theorem \ref{convert}, 

\[ 1+\sum_{\vec{a}} \sum_{\pi \in Mult(\vec{a})}
u^{P_1(\pi)} \prod_{i \geq 1} x_i^{a_i} =
\prod_{k=1}^{\infty} \frac{1}{1-\frac{x_k u}{1-\sum_{j>k} x_j}} \]
Differentiating with respect to $u$ and setting $u=1$ implies that the
sought expectation is

\begin{eqnarray*}
&& \frac{1}{{n \choose a_1,a_2,\cdots}} Coeff. \ of \prod x_i^{a_i} \
in \sum_k \frac{x_k}{1-\sum x_i} \frac{1}{1-\sum_{j \geq k} x_j}\\
&=& \frac{1}{{n \choose a_1,a_2,\cdots}} \sum_{k: a_k>0} \sum_{b_k,b_k+1,\cdots \geq 0}
 {a_1+\cdots+a_{k-1}+b_k+b_{k+1}+\cdots \choose a_1,\cdots,a_{k-1},b_k,b_{k+1},\cdots}\\
&& {a_k-1-b_k + a_{k+1}-b_{k+1} + \cdots \choose a_k-1-b_k, a_{k+1}-b_{k+1},\cdots}\\
&=& \frac{1}{{n \choose a_1,a_2,\cdots}} \sum_{k: a_k>0} \frac{1}{a_1!\cdots a_{k-1}! (a_k-1)! a_{k+1}! \cdots} \sum_{b_k,b_{k+1},\cdots \geq 0} \frac{{a_k-1 \choose b_k}{a_{k+1} \choose b_{k+1}}{a_{k+2} \choose b_{k+2}} \cdots}{{\sum a_i -1 \choose a_1+\cdots+a_{k-1}+b_k+b_{k+1}+\cdots}}
\end{eqnarray*}

	Letting $s=b_k+b_{k+1}+\cdots$, this simplifies to

\begin{eqnarray*}
&& \sum_{k: a_k>0} \sum_{s=0}^{a_k+a_{k+1}+\cdots-1} \frac{1}{\sum
a_i} \frac{{a_k+a_{k+1}+\cdots-1 \choose s}}{{\sum a_i -1 \choose
a_1+\cdots +a_{k-1}+s}}\\
&=& \sum_{k: a_k>0} \frac{1}{\sum a_i {\sum a_i -1 \choose
a_1+\cdots +a_{k-1}}} \sum_{s=0}^{a_k+a_{k+1}+\cdots-1} {a_1+\cdots+a_{k-1}+s \choose s}\\
&=& \sum_{k: a_k>0} \frac{a_k}{a_1+\cdots +a_{k-1}+1}\\
&=& \sum_k \frac{a_k}{a_1+\cdots +a_{k-1}+1}
\end{eqnarray*}

	The second calculation is similar.
\end{proof}

	As a final result, we study patience sorting applied to
$I_{2n}$, the fixed point free involutions in the symmetric group
$S_{2n}$. By \cite{R}, the number of piles in such a game relates to
the eigenvalues of random symplectic and orthogonal
matrices. Consequently this restricted version of patience sorting
merits further study. Proposition \ref{factorinv} shows that the
generating function for the first pile size factors.

\begin{prop} \label{factorinv}
 \[ \sum_{\pi \in I_{2n}} x^{P_1(\pi)} =
\prod_{i=1}^n (x^2+2(i-1)) \]
\end{prop}

\begin{proof} The proposition is proved by induction, the base case
being trivial. Suppose that the proposition holds for
$I_{2(n-1)}$. Given $\pi \in I_{2n}$ let $j$ be the symbol with which
$2n$ is switched. If $j \neq 1$, then $P_1(\pi)$ is the same as
$P_1(\pi')$ where $\pi'$ is obtained by crossing the symbols $j$ and
$2n$ out of $\pi$. If $j=1$, then $P_1(\pi)=P_1(\pi')+2$, where
$P_1(\pi')$ is obtained by crossing the symbols $1,2n$ out of
$\pi$. Consequently,

\[ \sum_{\pi \in I_{2n}} x^{P_1(\pi)} = (2n-2) \sum_{\pi \in I_{2n-2}}
x^{P_1(\pi)} + x^2 \sum_{\pi \in I_{2n-2}} x^{P_1(\pi)} \] and the
result follows by induction. \end{proof}

\section{Acknowledgements} This research was supported by an NSF
Postdoctoral Fellowship. The author thanks Persi Diaconis for ongoing
discussions.

\end{document}